\documentclass[twoside]{article} %% -*-latex-*-

\usepackage{jwn-article}
\usepackage{fullpage}
\usepackage{amsfonts}
\usepackage{mathptmx}
\usepackage{colortbl}
\arrayrulecolor[gray]{0.67}

\usepackage{hyperref}
\usepackage{graphicx}
\usepackage{moreverb}
\usepackage{alltt}
\usepackage[all]{xy}

\usepackage{amssymb}
\renewcommand{\ge}{\geqslant}
\renewcommand{\le}{\leqslant}

\newcommand{\cmd}[1]{\mathrm{#1}}
\newcommand{\prob}[2][]{\ensuremath{\cmd{Pr}_{#1}\left(#2\right)}}

\newcommand{\condprob}[3][]{\prob[#1]{#2\,|\,#3}}

\newcommand{\state}[1]{\ensuremath{\mathsf{#1}}}

\newcommand{\tval}{\ensuremath{\mathsf{T}}}
\newcommand{\fval}{\ensuremath{\mathsf{F}}}

\renewcommand{\phi}{\varphi}

\usepackage{stmaryrd}

\newcommand{\rname}[1]{\textsc{\MakeLowercase{#1}}}

\renewcommand{\theenumi}{\normalfont(\roman{enumi})}

\renewcommand{\theenumii}{\normalfont(\alph{enumii})}

\newcommand{\thd}[1]{\bfseries{#1}} % table heading

\title{Aristotle's Logic Computed by Parametric~Probability and
  Linear~Optimization}

\author{Joseph W. Norman \\ University of Michigan \\
\texttt{jwnorman@umich.edu}}
\date{June, 2014}

\begin{document}

\maketitle

\section{Introduction}

A new computational method is presented to implement the system of
deductive logic described by Aristotle in \emph{Prior Analytics}
\cite{aristotle-prior}.  Each Aristotelian problem is interpreted as a
parametric probability network in which the premises give constraints
on probabilities relating the problem's categorical terms (major,
minor, and middle).  Each probability expression from this network is
evaluated to yield a linear function of the parameters in the
probability model.  By this approach the constraints specified as
premises translate into linear equalities and inequalities involving a
few real-valued variables.  The problem's figure (schema) describes
which specific probabilities are constrained, relative to those that
are queried.  Using linear optimization methods, the minimum and
maximum feasible values of certain queried probabilities are computed,
subject to the constraints given as premises.  These computed
solutions determine precisely which conclusions are necessary
consequences of the premises.  In this way, Aristotle's logical
deductions can be accomplished by means of numerical computation.

This work is a synthesis of several existing methods, with the
addition of a few new ideas.  The most relevant prior work is that of
Boole, who presented several innovations in his 1854 treatise on the
\emph{Laws of Thought} \cite{boole}.  Boole demonstrated that logical
propositions can be represented as algebraic formulas; more
specifically, that statements of what we now call propositional
calculus can be expressed as polynomials with real-number
coefficients.  Boole showed how to compute interesting results about
logical propositions by solving systems of polynomial equations.
Boole also showed useful relationships between statements of logic and
statements of probability.  In the late 20th century Pearl and others
developed techniques for graphical probability models (Bayesian
networks) which offer several benefits regarding representation and
inference \cite{pearl,koller-friedman}.  Several investigators have
described methods for symbolic inference in probability networks;
these methods can be used to calculate polynomial formulas for queried
probability expressions
\cite{darwiche,castillo-networks,dambrosio,ideal-lisp}.  Boole already
formulated optimization problems with polynomial objectives and
constraints derived from probability expressions; he solved his
problems by \emph{ad hoc} algebraic manipulations \cite{boole}.  Today
there are general methods for solving linear and nonlinear polynomial
optimization problems.  For the linear case, efficient computational
methods were developed in the middle of the 20th century
\cite{dantzig,luenberger}.  It remains a challenge to compute exact
global solutions to unrestricted nonlinear polynomial optimization
problems; there are many promising methods which use various kinds of
approximation \cite{bertsekas,sherali,li,norman-nlp}.

Two new ideas are presented here which complement these existing
methods.  First, a taxonomy of Aristotelian categorical statements is
developed, with a distinction drawn between `primary' and `composite'
relations.  In this taxonomy, primary relations are mutually
exclusive, whereas composite relations may overlap.  For both kinds of
relations the case of an impossible subject (antecedent) is handled
explicitly.  This taxonomy results in several more types of
categorical statements than are usually considered (we end up with
seven).  One benefit is that existential fallacies are prevented.
Second, the concept of `complementary' syllogism is introduced, to
contrast with `classical' syllogism.  In a complementary syllogism the
subject of the deduced categorical statement is held to be false
instead of true.  Complementary syllogisms let us extract additional
information from categorical premises that would otherwise be lost to
analysis.  These new ideas (primary versus composite categorical
relations; classical versus complementary syllogism) are combined with
the existing methods mentioned above (Boole's mathematical logic;
probability networks and symbolic probability inference; linear
optimization) to provide the method of analysis presented here.  Let
us call this synthesis the `probability-optimization paradigm' for
framing Aristotle's logic.

%This `probability-optimization paradigm' for framing Aristotle's logic
%follows Boole's mathematical treatment of logic and probability
%\cite{boole}.

\subsection{The Structure of Aristotle's Logical Problems}
\label{sec:review}

The problems in Aristotle's \emph{Prior Analytics} involve three
categorical terms, called `major', `middle', and `minor', each of
which can be either true or false.  Let us use $A$ for the major term,
$B$ for the middle term, and $C$ for the minor term.  We abbreviate
truth as $\tval$ and falsity as $\fval$.  The major and minor terms
are also called the `extreme' terms, in contrast to the middle.  Each
Aristotelian problem consists of two premises and a query.  Each
premise has a subject and a predicate, each of which is one of the
three categorical terms.  Each premise also has a type that quantifies
the relationship between its subject and predicate terms. As described
in detail in Section~\ref{sec:translate}, we shall recognize five main
types of relationships: \rname{universal-affirmative},
\rname{universal-negative}, \rname{particular-affirmative},
\rname{particular-negative}, and \rname{particular-intermediate}.  The
\rname{universal} relationships are further subdivided into \rname{material}
and \rname{existential} subtypes; this expands our repertoire of
categorical relationships to seven types.  These seven relations,
which are not mutually exclusive, are composites built from four
primary relations which are mutually exclusive; details appear in
Section~\ref{sec:primary-composite}.

Aristotle imposed a few restrictions on how the various categorical
terms may be used in a problem's premises and query.  Both of the
problem's premises must use its middle term $B$.  One premise must use
the problem's major term $A$, and the other premise must use the
problem's minor term $C$ (hence these premises are called `major' and
`minor' themselves).  Within each premise either term may occupy the
position of subject or of predicate.  These restrictions allow four
possible figures for Aristotelian problems, as shown in
Table~\ref{tbl:figures} and discussed further in
Section~\ref{sec:figures}.  Regardless of which figure is used in a
problem, the main query is always the same: to find out what
relationship between the major and minor terms is required by the
given premises.  In this query the major term $A$ is used as predicate
and the minor term $C$ used as subject; this choice of positions is
precisely what distinguishes major from minor.  The inferred solution
is a subset of the seven types of categorical relationships introduced
above.  We shall say that a `syllogism' (deduction) has occurred when
at least one of these types of relationships must hold.  On the other
hand there is no syllogism when the premises do not require any
particular relationship between the major and minor terms.

Note that mathematical relationships may be asserted as constraints or
derived as solutions; these are two different roles.  For example, it
is one thing to assert the relation $x=2$ as a constraint that should
be satisfied, and a different thing to derive the relation $x=2$ as a
solution to some other system of constraints (for example, as one of
the two real solutions to the equation $x^2=4$).  Note also that
constraints are not commandments; the constraint $x=2$ does not
guarantee the solution $x=2$, for the complete system of equations
could be inconsistent with no solutions at all.  Anyway, we shall view
Aristotle's logical problems as systems of equations, both
philosophically and practically.  Philosophically, we shall regard
categorical statements like `$A$ belongs to some $B$' as relations
like $x=2$ or $y>0$, to be used in either of the two roles just
mentioned: sometimes asserted as constraints, and sometimes derived as
solutions to other constraints.  Practically, we shall translate
logical statements about true/false terms into algebraic equalities
and inequalities involving real-valued variables (through the
intermediate device of a probability model), and then use standard
algebraic and numerical methods to compute solutions to the original
logical problems.

\subsection{Notation for Problems in Four Figures}
\label{sec:figures}

\begin{table}
  \sf
  \begin{tabular}{l|ccccc} \hline
    & \bfseries Major & \bfseries Minor & \\
    \bfseries Figure & 
    \bfseries Premise & 
    \bfseries Premise &
    \bfseries Query &
    \bfseries Flat Diagram &
    \bfseries Triangular Diagram 
    \\ \hline\hline

    First & $A\mathfrak{m}B$ & $B\mathfrak{n}C$ & $A\mathfrak{s}C$ &
    \begin{math}
      \vcenter{\xymatrix{
        A & 
        B \ar@{>}[l]_{\mathfrak{m}} & 
        C \ar@{>}[l]_{\mathfrak{n}} \ar@{.>}@/^2ex/[ll]^{\mathfrak{s}}
      }}
    \end{math}
    &
    \begin{math}
      \vcenter{\xymatrix @R=1.5em @C=0.5em {
          & 
          B \ar@{>}[dl]_{\mathfrak{m}} 
          & \\
          A
          & & 
          C \ar@{.>}[ll]^{\mathfrak{s}} 
          \ar@{>}[ul]_{\mathfrak{n}}
      }}
    \end{math}
    \\ \hline

    Second & $B\mathfrak{m}A$ & $B\mathfrak{n}C$ & $A\mathfrak{s}C$ &
    \begin{math}
      \vcenter{\xymatrix{
        B & 
        A \ar@{>}[l]_{\mathfrak{m}} & 
        C \ar@{>}@/^2ex/[ll]^{\mathfrak{n}} \ar@{.>}[l]_{\mathfrak{s}}
      }}
    \end{math}
    &
    \begin{math}
      \vcenter{\xymatrix @R=1.5em @C=0.5em {
          & 
          B
          & \\
          A
          \ar@{>}[ur]^{\mathfrak{m}} 
          & & 
          C \ar@{.>}[ll]^{\mathfrak{s}} 
          \ar@{>}[ul]_{\mathfrak{n}}
      }}
    \end{math}
    \\ \hline

    Third & $A\mathfrak{m}B$ & $C\mathfrak{n}B$ & $A\mathfrak{s}C$ &
    \begin{math}
      \vcenter{\xymatrix{
        A & 
        C \ar@{.>}[l]_{\mathfrak{s}} & 
        B \ar@{>}[l]_{\mathfrak{n}} \ar@{>}@/^2ex/[ll]^{\mathfrak{m}}
      }}
    \end{math}
    &
    \begin{math}
      \vcenter{\xymatrix @R=1.5em @C=0.5em {
          & 
          B
          \ar@{>}[dl]_{\mathfrak{m}} 
          \ar@{>}[dr]^{\mathfrak{n}}
          & \\
          A
          & & 
          C \ar@{.>}[ll]^{\mathfrak{s}} 
      }}
    \end{math}
    \\ \hline

    Fourth & $B\mathfrak{m}A$ & $C\mathfrak{n}B$ & $A\mathfrak{s}C$ &
    \begin{math}
      \vcenter{\xymatrix{
        C \ar@{.>}@/_2ex/[rr]_{\mathfrak{s}} &
        B \ar@{>}[l]_{\mathfrak{n}} & 
        A \ar@{>}[l]_{\mathfrak{m}}
      }}
    \end{math}
    &
    \begin{math}
      \vcenter{\xymatrix @R=1.5em @C=0.5em {
          & 
          B
          \ar@{>}[dr]^{\mathfrak{n}}
          & \\
          A \ar@{>}[ur]^{\mathfrak{m}} 
          & & 
          C \ar@{.>}[ll]^{\mathfrak{s}} 
      }}
    \end{math}
    \\ \hline
  \end{tabular}
  \caption{Four figures of Aristotelian problems, using major term
    $A$, middle term $B$, and minor term $C$.  The major premise has
    type $\mathfrak{m}$, the minor premise has type $\mathfrak{n}$,
    and the goal is to find the implied types $\mathfrak{s}$ for the
    query statement $A\mathfrak{s}C$.  The possible relation types are
    defined in Table~\ref{tbl:composites}.  In the diagrams each arrow
    points from the subject to the predicate of a conditional
    statement; solid arrows indicate premises and dotted arrows
    indicate queries.}
  \label{tbl:figures}
\end{table}

%Before discussing the details of how probability and optimization
%results indicate logical deductions, let 
We now consider textual and graphical notation for Aristotle's logical
problems.  Table~\ref{tbl:figures} describes the four figures of
Aristotelian problems in symbolic and graphical notation.  The
symbolic notation indicates which term is the predicate and which is
the subject of each categorical statement (the two statements asserted
as premises, and the one statement used as a query).  Capital italic
letters $A$, $B$, and $C$ stand for the major, middle, and minor
terms.  Small italic letters (with optional accent marks) from the set
\begin{math}
  \{ a, \acute{a}, e, \acute{e}, i, o, u \}
\end{math}
stand for types of categorical relationships, which are defined by
constraints on probabilities as shown in Table~\ref{tbl:composites}.
For example $a$ stands for the \rname{universal-affirmative-material}
relation and $u$ stands for the \rname{particular-intermediate}
relation.  Gothic letters $\mathfrak{m}$, $\mathfrak{n}$, and
$\mathfrak{s}$ stand for categorical relations from the set
\begin{math}
  \{ a, \acute{a}, e, \acute{e}, i, o, u \}
\end{math}: $\mathfrak{m}$ for the
relation type of the major premise, $\mathfrak{n}$ for the relation
type of the minor premise, and $\mathfrak{s}$ for the relation type of
the queried statement.  For example, a problem in the second figure
has a major premise $B\mathfrak{m}A$ with predicate $B$, subject $A$,
and relation type $\mathfrak{m}$; it has a minor premise
$B\mathfrak{n}C$ with predicate $B$, subject $C$, and relation type
$\mathfrak{n}$; and it has query $A\mathfrak{s}C$ with predicate $A$,
subject $C$, and relation type $\mathfrak{s}$.  The example problem
from Section~\ref{sec:linear} follows the second figure.  Its major
premise $BeA$ uses the \rname{Universal-negative-material} relation
(denoted $e$) as $\mathfrak{m}$, and the minor premise $BiC$ uses the
\rname{Particular-affirmative} relation (denoted $i$) as
$\mathfrak{n}$.  The meaning of the query $A\mathfrak{s}C$ is
discussed in Section~\ref{sec:criteria}.

It is interesting that Aristotle already described his figures in
graphical language, indicating the positions of the various
categorical terms on the page.  For example, regarding his second
figure, Aristotle wrote:
\begin{quote}
  \ldots by middle term in it I mean that which is predicated by both
  subjects, by extremes the terms of which this is said, by major
  extreme that which lies near the middle, by minor that which is
  further away from the middle.  The middle term stands outside the
  extremes, and is first in position. (\cite{aristotle-prior}~26b35)
\end{quote}
Furthermore, Aristotle described his third figure in this way:
\begin{quote}
  \ldots by extremes I mean the predicates, by the major extreme that
  which is further from the middle, by the minor that which is nearer
  to it.  The middle term stands outside the extremes, and is last in
  position. (\cite{aristotle-prior}~28a15)
\end{quote}
The graphical diagrams included in Table~\ref{tbl:figures} realize
Aristotle's original textual descriptions in one view (the `flat'
diagrams), and use a different graph layout in an alternative view
(the `triangular' diagrams).  As has become customary, a fourth figure
has been added.  The fourth figure is related to the first by having
the major and minor terms swapped with one another.

\subsection{A Basic Probability Model}
\label{sec:basic-model}

\begin{table}
%  \small\sf
\begin{tabular}[c]{lll|l}\hline
\multicolumn{1}{l}{$A$} & 
\multicolumn{1}{l}{$B$} & 
\multicolumn{1}{l|}{$C$} & 
\multicolumn{1}{l}{$\prob[0]{{A, B, C}}$} \\ \hline\hline
$\state{T}$ & 
$\state{T}$ & 
$\state{T}$ & 
$x_{1}$ \\ \hline
$\state{T}$ & 
$\state{T}$ & 
$\state{F}$ & 
$x_{2}$ \\ \hline
$\state{T}$ & 
$\state{F}$ & 
$\state{T}$ & 
$x_{3}$ \\ \hline
$\state{T}$ & 
$\state{F}$ & 
$\state{F}$ & 
$x_{4}$ \\ \hline
$\state{F}$ & 
$\state{T}$ & 
$\state{T}$ & 
$x_{5}$ \\ \hline
$\state{F}$ & 
$\state{T}$ & 
$\state{F}$ & 
$x_{6}$ \\ \hline
$\state{F}$ & 
$\state{F}$ & 
$\state{T}$ & 
$x_{7}$ \\ \hline
$\state{F}$ & 
$\state{F}$ & 
$\state{F}$ & 
$x_{8}$ \\ \hline
\end{tabular}
  \caption{The input probability table $\prob[0]{A,B,C}$ for the basic
    model, with real-valued parameters $x_1$ through $x_8$ subject to
    the constraints $0 \le x_i \le 1$ and $\sum_i x_i = 1$.}
  \label{tbl:basic-cpt}
\end{table}

Our basic probability model represents the joint probabilities of the
three categorical terms $A$, $B$, and $C$.  With two possible truth
values for each of the three terms, there are $2^3$ or $8$ possible
combinations of truth values.  To each combination $i$ of truth values
we assign a symbolic parameter $x_i$ that represents its probability.
These parameters and their associated combinations of truth values are
shown in Table~\ref{tbl:basic-cpt} as the input probability table
$\prob[0]{A,B,C}$.  To respect the laws of probability, these
parameters
\begin{math}
x_{1}
,\ldots,
x_{8}
\end{math}
are constrained by $0 \le x_i \le 1$ and $\sum_i x_i = 1$.  Parametric
probability networks such as this basic model are used like databases
to answer queries.  Each query requests an unconditioned probability
or a conditional probability.  Each response is a polynomial or a
quotient of polynomials in the model's parameters.  For example,
starting from the inputs shown in Table~\ref{tbl:basic-cpt}, the
respective probabilities that $B$ is true, that $B$ and $A$ are both
true, that $B$ is true and $A$ is false, and that $A$ is true given
that $B$ is true are computed as the following algebraic expressions:
\begin{eqnarray}
  \prob{B} & \Rightarrow &
x_{1} + x_{2} + x_{5} + x_{6}
  \\
  \prob{B,A} & \Rightarrow &
x_{1} + x_{2}
  \\
  \prob{B,\overline{A}} & \Rightarrow &
x_{5} + x_{6}
  \\
  \condprob{A}{B} & \Rightarrow &
\left(x_{1} + x_{2}\right) / \left(x_{1} + x_{2} + x_{5} + x_{6}\right)
\end{eqnarray}
As you can see, each of these calculated values is either the sum of
several input probabilities from Table~\ref{tbl:basic-cpt} or the
quotient of two such sums.  Table~\ref{tbl:basic-output} shows several
output probabilities computed from the inputs in
Table~\ref{tbl:basic-cpt}.  These outputs will be useful for the
analysis that follows.  It happens with this basic probability model
that all computed probabilities are \emph{linear} functions of the
$x_i$ parameters (or quotients of such linear functions).  Other
probability models can yield \emph{nonlinear} polynomials and
quotients (when the full-joint probability has been factored into
multiple input tables).

\begin{table}
  \sf
  \begin{tabular}{rl@{\qquad}rl@{\qquad}rl}
    (a) &
\begin{tabular}[c]{ll|l}\hline
\multicolumn{1}{l}{$A$} & 
\multicolumn{1}{l|}{$B$} & 
\multicolumn{1}{l}{$\prob{{A, B}}$} \\ \hline\hline
$\state{T}$ & 
$\state{T}$ & 
$x_{1} + x_{2}$ \\ \hline
$\state{T}$ & 
$\state{F}$ & 
$x_{3} + x_{4}$ \\ \hline
$\state{F}$ & 
$\state{T}$ & 
$x_{5} + x_{6}$ \\ \hline
$\state{F}$ & 
$\state{F}$ & 
$x_{7} + x_{8}$ \\ \hline
\end{tabular}
    & (b) &
\begin{tabular}[c]{ll|l}\hline
\multicolumn{1}{l}{$B$} & 
\multicolumn{1}{l|}{$C$} & 
\multicolumn{1}{l}{$\prob{{B, C}}$} \\ \hline\hline
$\state{T}$ & 
$\state{T}$ & 
$x_{1} + x_{5}$ \\ \hline
$\state{T}$ & 
$\state{F}$ & 
$x_{2} + x_{6}$ \\ \hline
$\state{F}$ & 
$\state{T}$ & 
$x_{3} + x_{7}$ \\ \hline
$\state{F}$ & 
$\state{F}$ & 
$x_{4} + x_{8}$ \\ \hline
\end{tabular}
    & (c) &
\begin{tabular}[c]{ll|l}\hline
\multicolumn{1}{l}{$C$} & 
\multicolumn{1}{l|}{$A$} & 
\multicolumn{1}{l}{$\prob{{C, A}}$} \\ \hline\hline
$\state{T}$ & 
$\state{T}$ & 
$x_{1} + x_{3}$ \\ \hline
$\state{T}$ & 
$\state{F}$ & 
$x_{5} + x_{7}$ \\ \hline
$\state{F}$ & 
$\state{T}$ & 
$x_{2} + x_{4}$ \\ \hline
$\state{F}$ & 
$\state{F}$ & 
$x_{6} + x_{8}$ \\ \hline
\end{tabular}
  \\ \\
  (d) &
\begin{tabular}[c]{l|l}\hline
\multicolumn{1}{l|}{$A$} & 
\multicolumn{1}{l}{$\prob{{A}}$} \\ \hline\hline
$\state{T}$ & 
$x_{1} + x_{2} + x_{3} + x_{4}$ \\ \hline
$\state{F}$ & 
$x_{5} + x_{6} + x_{7} + x_{8}$ \\ \hline
\end{tabular}
  & (e) &
\begin{tabular}[c]{l|l}\hline
\multicolumn{1}{l|}{$B$} & 
\multicolumn{1}{l}{$\prob{{B}}$} \\ \hline\hline
$\state{T}$ & 
$x_{1} + x_{2} + x_{5} + x_{6}$ \\ \hline
$\state{F}$ & 
$x_{3} + x_{4} + x_{7} + x_{8}$ \\ \hline
\end{tabular}
  & (f) &
\begin{tabular}[c]{l|l}\hline
\multicolumn{1}{l|}{$C$} & 
\multicolumn{1}{l}{$\prob{{C}}$} \\ \hline\hline
$\state{T}$ & 
$x_{1} + x_{3} + x_{5} + x_{7}$ \\ \hline
$\state{F}$ & 
$x_{2} + x_{4} + x_{6} + x_{8}$ \\ \hline
\end{tabular}
  \end{tabular}
  \caption{A few output probability tables computed from the inputs in
    Table~\ref{tbl:basic-cpt}.}
  \label{tbl:basic-output}
\end{table}

The essential methods of symbolic probability inference were described
well enough several centuries ago \cite{demoivre}.  There have since
been developed more rigorous mathematical formulations, more efficient
inference algorithms, and powerful graphical models
\cite{kolmogorov,jensen,koller-friedman}.  The author has developed
some computational methods for parametric probability networks as well
\cite{norman-thesis,norman-problogic}; these methods include some
idiosyncratic notation that is reviewed presently.  \emph{Input} and
\emph{output} probabilities are distinguished from one another.  Input
probabilities, used to specify the probability model, are written with
the subscript $0$, as in $\prob[0]{A,B,C}$.  Output probabilities,
computed from the inputs, are written with no subscript, as in
$\prob{B}$.  The double right arrow $\Rightarrow$ is used to indicate
computation, such as the evaluation of a symbolic probability
expression or the simplification of an arithmetical formula.  This
meaning is distinct from the test or assertion of equality denoted
with the usual equal sign $=$.  Finally, probability tables and their
elements share similar notation.  A probability expression such as
$\prob{A,B}$ may refer to a table containing several values, such as
the four elements shown as Table~\ref{tbl:basic-output}~(a).  But we
can also use for example $A$ to abbreviate the event $A=\tval$ and
$\overline{A}$ to abbreviate the event $A=\fval$, and hence use
$\prob{A,B}$ to mean the individual element $\prob{A=\tval,B=\tval}$.
The default used here is that probability expressions refer to
individual elements; it will be announced in the neighboring text when
a probability expression refers instead to an entire table containing
several elements.

\section{From Categorical Statements to Linear Equalities and
  Inequalities} 
\label{sec:translate}

\subsection{Naive Types of Categorical Relationships}
\label{sec:naive}

Let us now translate Aristotelian categorical statements into linear
equalities and inequalities involving the parameters of the basic
probability model from Section~\ref{sec:basic-model}.  To begin, we
regard a categorical statement with predicate $P$ and subject $Q$ as a
relation involving $\condprob{P}{Q}$, the conditional probability that
$P$ is true given that $Q$ is true.  Here $P$ and $Q$ can be any of
the three categorical terms $A$, $B$, or $C$.  At first glance,
Aristotle's \emph{Prior Analytics} describes four types of relations
between categorical terms, which correspond the listed
conditional-probability statements:
\begin{equation}
  \begin{tabular}{lll}
    \rname{Universal-affirmative} &
    `$P$ belongs to all $Q$' &
    $\condprob{P}{Q}=1$
    \\
    \rname{Universal-negative} &
    `$P$ belongs to no $Q$' &
    $\condprob{P}{Q}=0$
    \\
    \rname{Particular-affirmative} &
    `$P$ belongs to some $Q$' &
    $\condprob{P}{Q}>0$
    \\
    \rname{Particular-negative} &
    `$P$ does not belong to some $Q$' &
    $\condprob{P}{Q}<1$
  \end{tabular}
  \label{eq:pre-primaries}
\end{equation}
Recall that conditional probabilities are defined as quotients of
unconditioned probabilities:
\begin{eqnarray}
  \condprob{P}{Q} & \equiv & 
  \frac{\prob{P,Q}}{\prob{Q}}
  \label{eq:condprob}
\end{eqnarray}
There are two troublesome issues with the four types of relations
listed above.  The first issue is that there is no prescription for
how to handle the case that $\prob{Q}=0$ (meaning that it is
impossible \emph{a priori} for the subject term $Q$ to be true).
Since the laws of probability require that $\prob{P,Q}=0$ when
$\prob{Q}=0$, this exceptional case would force the quotient shown in
Equation~\ref{eq:condprob} to have the indefinite value $0/0$.  It is
ambiguous whether equations such as $0/0=0$ and $0/0=1$ should be
considered satisfied (on the one hand, both rearranged equations $0 =
0 \cdot 0$ and $0 = 0 \cdot 1$ are true; on the other hand, neither
$0$ nor $1$ provides a \emph{unique} solution to the rearranged
equation $0=0 \cdot c$).  The second issue is that the listed
relations are not mutually exclusive.  For example $\condprob{P}{Q}=1$
requires also $\condprob{P}{Q}>0$, and conversely $\condprob{P}{Q}>0$
leaves it possible but not certain that $\condprob{P}{Q}=1$.  This
lack of exclusivity may lead to confusion about precisely which
relations hold true in any given circumstance.

As an aside, note some potential confusion regarding the negation of a
\rname{particular-affirmative} premise.  The English phrase `$P$ does
not belong to some $Q$' leaves ambiguity about what is negated.  This
phrase could be interpreted to mean, `It is not the case that $P$
belongs to some $Q$'---suggesting the conditional-probability
constraint $\condprob{P}{Q}=0$.  Or this phrase could be interpreted
to mean, `The negation of $P$ belongs to some $Q$'---suggesting the
constraint $\condprob{\overline{P}}{Q}>0$ (or its equivalent
$\condprob{P}{Q}<1$ which is listed above).  We shall assume the
latter of these interpretations.

\subsection{Primary and Composite Relations}
\label{sec:primary-composite}

In order to address the troublesome issues with the naive categorical
relations presented above, let us develop a different initial set of
relation types---making direct use of the numerator $\prob{P,Q}$ and
denominator $\prob{Q}$ from Equation~\ref{eq:condprob} instead of
their quotient $\condprob{P}{Q}$.  Let us identify these four
\emph{primary relations} concerning a predicate term $P$ and a subject
term $Q$:
\begin{equation}
  \begin{tabular}{rll}
    R1. & \rname{Impossible-subject} &
    `There are no $Q$'
    \\
    R2. & \rname{Universal-negative-existential} &
    `$P$ belongs to no $Q$, and there are some $Q$'
    \\
    R3. & \rname{Particular-intermediate} &
    `$P$ belongs to some but not all $Q$'
    \\
    R4. & \rname{Universal-affirmative-existential} &
    `$P$ belongs to all $Q$, and there are some $Q$'
  \end{tabular}
  \label{eq:primaries}
\end{equation}
We shall define these primary relations using constraints on two
probabilities from the model in Section~\ref{sec:basic-model}: the
probability $\prob{Q}$ that the subject term $Q$ is true, and the
probability $\prob{P,Q}$ that both the predicate term $P$ and the
subject term $Q$ are true.  We consider the two cases $\prob{Q}=0$ and
$\prob{Q}>0$ on one axis, and the three cases $\prob{P,Q}=0$, $0 <
\prob{P,Q} < \prob{Q}$, and $\prob{P,Q}=\prob{Q}$ on the other.  This
gives the following matrix of constraints, whose entries define the
primary relations described in Equation~\ref{eq:primaries}:
\begin{equation}
  \begin{array}{|r||c|c|c|} \hline
    & \prob{P,Q}=0 & 0 < \prob{P,Q} < \prob{Q} & \prob{P,Q}=\prob{Q} 
    \\ \hline\hline
    \prob{Q}=0 & \mathrm{R1} & & \mathrm{R1} \\ \hline
    \prob{Q}>0 & \mathrm{R2} & \mathrm{R3} & \mathrm{R4} \\ \hline
  \end{array}
  \label{eq:primary-matrix}
\end{equation}
For example R1 is defined as the case that $\prob{Q}=0$; and R2 is
defined as the case that $\prob{Q}>0$ and $\prob{P,Q}=0$.  Note that
when $\prob{Q}=0$ the laws of probability require $\prob{P,Q}=0$ also.
Hence in the case R1 both constraints $\prob{P,Q}=0$ and
$\prob{P,Q}=\prob{Q}$ are satisfied.  Furthermore, when $\prob{Q}=0$
it is not possible for $\prob{P,Q}$ to be strictly greater nor
strictly less than zero; hence the corresponding entry of the matrix
in Equation~\ref{eq:primary-matrix} is empty.  The constraints of
Equation~\ref{eq:primary-matrix} imply derived constraints on the
conditional probability $\condprob{P}{Q}$: R1 requires that
$\condprob{P}{Q}$ must have the indefinite value $0/0$; R2 requires
that $\condprob{P}{Q}$ must equal $0$; R4 requires that
$\condprob{P}{Q}$ must equal $1$; and R3 requires that
$\condprob{P}{Q}$ must be strictly greater than $0$ and strictly less
than $1$.

The primary relations R1, R2, R3, and R4 given by Equations
\ref{eq:primaries} and \ref{eq:primary-matrix} address both issues
identified in Section~\ref{sec:naive}: the defining constraints are
mutually exclusive and they include the case $\prob{Q}=0$ explicitly.
We can now proceed to use various subsets of these four primary
relations to define \emph{composite relations} which include the seven
types of categorical relationships promised in the introduction.
These compositions are logical disjunctions.  For example, the
combination `R1 or R4' yields the composite relation that either
$\prob{Q}=0$, or $\prob{Q}>0$ and $\prob{P,Q}=\prob{Q}$: in other
words, the combined statement that `Either there are no $Q$, or there
are some $Q$ and $P$ belongs to all of them'.  This composite relation
can be specified as the integrated constraint $\prob{P,Q}=\prob{Q}$.
Recalling Equation~\ref{eq:condprob}, it follows from this integrated
constraint that the conditional probability $\condprob{P}{Q}$ must
have either the definite value $1$ or the indefinite value $0/0$
(depending on whether $\prob{Q}>0$ or $\prob{Q}=0$).  In a sense this
composite relation parallels the statement of material implication $Q
\rightarrow P$ from the propositional calculus; thus we call it the
\rname{Universal-affirmative-material} relation.

Table~\ref{tbl:composites} shows seven composite relations defined as
disjunctive combinations of the primary relations introduced above.
This set of composite relations is meant to be expressive rather than
exhaustive; in total there are $2^4$ or $16$ possible sets of $4$
primary relations.  Table~\ref{tbl:descriptions} gives
conditional-probability and natural-language descriptions of the seven
selected composite relations.  In both tables, the composite relations
are assigned codes to abbreviate them, based on the letters introduced
in medieval times to designate different types of Aristotelian
premises.  The traditional codes $a$, $e$, $i$, and $o$ are
supplemented with accented characters $\acute{a}$ and $\acute{e}$ that
distinguish \rname{existential} subtypes of \rname{universal}
statements from their \rname{material} counterparts.  Also, the letter
$u$ has been added to designate the \rname{particular} relation
meaning `some but not all'.  Within this document, for abbreviations
$PaQ$, $PiQ$, and so on, the predicate $P$ is displayed before the
relation code, and the subject $Q$ after it.  Be aware that some
authors use the opposite convention.  Aristotle's original texts did
not use such abbreviations at all.

\begin{table}
  \sf
  \begin{tabular}{l|c||c|c|c|c||l} \hline
    & & \multicolumn{4}{|c|}{\thd{Primary Relations}} &
    \\ \cline{3-6}
    \thd{Composite Relation} &
    \thd{Code} &
    \thd{R1} &
    \thd{R2} &
    \thd{R3} &
    \thd{R4} & 
    \thd{Probability Constraint} \\ \hline\hline
    
    Universal-affirmative-material & $PaQ$ & 
    $\bullet$ & & & $\bullet$ &
    $\prob{P,Q}=\prob{Q}$
    \\ \hline

    Universal-affirmative-existential & $P\acute{a}Q$ & 
    & & & $\bullet$ &
    $\prob{P,Q}=\prob{Q}$ and $\prob{Q}>0$
    \\ \hline

    Universal-negative-material & $PeQ$ &
    $\bullet$ & $\bullet$ & & &
    $\prob{P,Q}=0$
    \\ \hline

    Universal-negative-existential & $P\acute{e}Q$ &
    & $\bullet$ & & &
    $\prob{P,Q}=0$ and $\prob{Q}>0$
    \\ \hline

    Particular-affirmative & $PiQ$ &
    & & $\bullet$ & $\bullet$ & 
    $\prob{P,Q}>0$
    \\ \hline

    Particular-negative & $PoQ$ &
    & $\bullet$ & $\bullet$ & &
    $\prob{P,Q}<\prob{Q}$
    \\ \hline

    Particular-intermediate & $PuQ$ &
    & & $\bullet$ & &
    $\prob{P,Q}>0$ and $\prob{P,Q}<\prob{Q}$
    \\ \hline

%%     Impossible-subject & $\phantom{P}yQ$
%%     & $\bullet$ & & & & 
%%     $\prob{Q}=0$
%%     \\ \hline
    
  \end{tabular}
  \caption{Composite relations between a categorical predicate $P$ and
    a subject $Q$, using disjunctions of the primary relations R1, R2,
    R3, and R4 from Equation~\ref{eq:primary-matrix}.  Each bullet
    $\bullet$ indicates that the given primary relation is included in
    the given composite: for example
    \rname{Universal-affirmative-material} (code $a$) holds if either
    primary relation R1 or R4 holds.  The integrated probability
    constraint that defines each composite categorical relation is
    shown.}
  \label{tbl:composites}
\end{table}

\begin{table}
  \sf
  \begin{tabular}{l|l|l|l} \hline
    \bfseries Composite Relation &
    \bfseries Code &
    \bfseries Derived Cond.\ Prob. &
    \bfseries Natural-Language Description
    \\ \hline\hline
    
    Universal-affirmative-material & $PaQ$ &
    $\condprob{P}{Q}=1 \; \mbox{or} \; 0/0$ &
    $P$ belongs to all $Q$, or there are no $Q$
    \\ \hline

    Universal-affirmative-existential & $P\acute{a}Q$ &
    $\condprob{P}{Q}=1$ &
    $P$ belongs to all $Q$, and there are some $Q$
    \\ \hline

    Universal-negative-material & $PeQ$  &
    $\condprob{P}{Q}=0 \; \mbox{or} \; 0/0$ &
    $P$ belongs to no $Q$, or there are no $Q$
    \\ \hline

    Universal-negative-existential & $P\acute{e}Q$ &
    $\condprob{P}{Q}=0$ &
    $P$ belongs to no $Q$, and there are some $Q$
    \\ \hline

    Particular-affirmative & $PiQ$ &
    $\condprob{P}{Q}>0$ &
    $P$ belongs to some $Q$
    \\ \hline

    Particular-negative & $PoQ$ &
    $\condprob{P}{Q}<1$ &
    The negation of $P$ belongs to some $Q$
    \\ \hline

    Particular-intermediate & $PuQ$ &
    $0<\condprob{P}{Q}<1$ &
    $P$ belongs to some but not all $Q$
    \\ \hline

%%     Impossible-subject & $\phantom{P}yQ$ &
%%     $\condprob{P}{Q}=0/0$ &
%%     There are no $Q$
%%     \\ \hline
    
  \end{tabular}
  \caption{Composite relations between a categorical predicate $P$ and
    a subject $Q$, described in terms of natural language and
    conditional-probability constraints derived from the
    unconditioned-probability constraints in Table~\ref{tbl:composites}.}
  \label{tbl:descriptions}
\end{table}

\subsection{Instantiation as Linear Equalities and Inequalities}
\label{sec:instantiation}

It remains to instantiate the composite relations defined in
Table~\ref{tbl:composites} into specific equalities and inequalities
involving the parameters
\begin{math}
x_{1}
,\ldots,
x_{8}
\end{math}
of the basic probability model from Section~\ref{sec:basic-model},
when particular categorical terms ($A$, $B$, or $C$) have been chosen
as the predicate $P$ and as the subject $Q$.  This instantiation is
accomplished by using the results of symbolic probability inference
shown in Table~\ref{tbl:basic-output} to supply algebraic formulas for
the relevant probability expressions.

For example, let us consider the
\rname{Universal-affirmative-existential} statement with predicate
term $A$ and subject term $B$.  As Table~\ref{tbl:descriptions} shows,
this categorical statement $A\acute{a}B$ says that `$A$ belongs to all
$B$, and there are some $B$'.  Using Table~\ref{tbl:composites} and
substituting $A$ for the predicate term $P$ and $B$ for the subject
term $Q$, this statement $A\acute{a}B$ is defined by the probability
relations $\prob{B,A}=\prob{B}$ and $\prob{B}>0$.
Table~\ref{tbl:basic-output} gives the algebraic formulas for the
relevant probability expressions.  The probability $\prob{B,A}$,
meaning $\prob{B=\tval,A=\tval}$, is given in the first row of the
computed table $\prob{A,B}$ which appears as part~(a) of
Table~\ref{tbl:basic-output}:
\begin{math}
  \prob{B,A} \Rightarrow
x_{1} + x_{2}
\end{math}.
The probability $\prob{B}$, meaning
$\prob{B=\tval}$, is given in the first row of the computed table
$\prob{B}$ which appears as part~(e) of Table~\ref{tbl:basic-output}:
\begin{math}
  \prob{B} \Rightarrow
x_{1} + x_{2} + x_{5} + x_{6}
\end{math}.
Substituting these values into the probability relations
$\prob{B,A}=\prob{B}$ and $\prob{B}>0$ derived from
Table~\ref{tbl:composites}, the categorical statement $A\acute{a}B$ is
therefore defined by the following algebraic relations:
\begin{eqnarray}
x_{1} + x_{2}
  & = &
x_{1} + x_{2} + x_{5} + x_{6}
  \label{eq:aab1}
  \\
x_{1} + x_{2} + x_{5} + x_{6}
  & > & 0
  \label{eq:aab2}
\end{eqnarray}
Equation~\ref{eq:aab1} simplifies to
\begin{math}
x_{5} + x_{6}
  = 0
\end{math}.
Substituting this result, Equation~\ref{eq:aab2} then simplifies to
\begin{math}
x_{1} + x_{2}
  > 0
\end{math}.
As Table~\ref{tbl:descriptions} says, the constraints
$\prob{B,A}=\prob{B}$ and $\prob{B}>0$ that define the categorical
statement $A\acute{a}B$ require that the conditional probability
$\condprob{A}{B}$ must have the definite value $1$.  This is evident
from inspecting Equations \ref{eq:aab1} and \ref{eq:aab2} along with
the algebraic formula computed for the conditional probability:
\begin{eqnarray}
  \condprob{A=\tval}{B=\tval} & \Rightarrow &
\left(x_{1} + x_{2}\right) / \left(x_{1} + x_{2} + x_{5} + x_{6}\right)
\end{eqnarray}
Substituting
\begin{math}
x_{5} + x_{6}
  = 0
\end{math}
from the simplified Equation~\ref{eq:aab1} and using
\begin{math}
x_{1} + x_{2}
  > 0
\end{math}
from the simplified Equation~\ref{eq:aab2}, it follows that this
quotient expressing $\condprob{A}{B}$ must have the definite value
$1$.

\section{Computing Classical and Complementary Syllogisms}
\label{sec:compute}

The methods of Section~\ref{sec:translate} enable the translation of
Aristotelian categorical statements into linear equalities and
inequalities involving the real-valued parameters
\begin{math}
x_{1}
,\ldots,
x_{8}
\end{math}
of the basic probability model from Section~\ref{sec:basic-model}.
Such equalities and inequalities can be asserted as constraints
themselves, or inferred as solutions to other systems of constraints.
We now turn to the task of performing such inference: using linear
optimization to compute some linear equalities and inequalities from
others.  These computed algebraic results can be used to deduce
logical conclusions from Aristotelian premises.  We shall define two
kinds of deductions: \emph{classical} syllogisms, which follow
Aristotle's original practice; and \emph{complementary} syllogisms,
which offer a useful variation on the theme.  As explained further in
Section~\ref{sec:criteria}, the difference is that for classical
syllogisms the minor term is held to be true, whereas for
complementary syllogisms the minor term is held to be false (in the
relevant queries, which ask how the major term must be predicated upon
the minor term).

\subsection{Linear Programming to Bound Queried Probabilities}
\label{sec:linear}

Consider an objective function $f$ and some constraint functions
$g_1,g_2,\ldots,g_m$, each of which is a linear function of a list
$\mathbf{x}=(x_1,\ldots,x_n$) of variables.  Standard linear
optimization methods can compute the minimum and maximum feasible
values of the objective $f(\mathbf{x})$ subject to linear equality and
inequality constraints such as $g_1(\mathbf{x})=0$, $g_2(\mathbf{x})
\ge 0$, and so on \cite{luenberger}.  Standard linear optimization
methods can also detect inconsistent constraints, a potential
exception that is important to recognize.\footnote{It happens that
  following Aristotle's restrictions on terms and premises, there are
  no infeasible problems.  However with more than two premises, or
  with terms that do not follow the typical major/minor/middle
  arrangement, it is quite possible to give inconsistent categorical
  statements.  One simple example is the premise $AoA$ with the
  \rname{particular-negative} relation type and the same term used as
  both subject and predicate.  This premise translates to the
  unsatisfiable probability constraint $\prob{A,A}<\prob{A}$ which
  says that the probability that $A$ is true is strictly less than
  itself.}

For probabilistic translations of Aristotelian problems it is
necessary to reason with strict inequalities ($<$ and $>$) as well as
weak ones ($\le$ and $\ge$).  However common linear optimization
methods treat all inequalities as weak.  The following
`epsilon-inequality reformulation' works around this limitation.  We
choose a constant value $\epsilon$ such as
\begin{math}
0.01
\end{math}
or $1 \times 10^{-6}$ that is small relative to $1$, and replace each
strict-inequality constraint $g(\mathbf{x}) > h(\mathbf{x})$ with a
weak-inequality constraint $g(\mathbf{x}) \ge h(\mathbf{x}) +
\epsilon$.  The optimization problem thus reformulated is solved using
standard linear programming methods to find the minimum and maximum
feasible values of the objective function $f(\mathbf{x})$ subject to
the given constraints.  Any computed minimum solution greater than or
equal to $\epsilon$ is interpreted to mean that the objective
$f(\mathbf{x})$ must be strictly greater than zero; likewise any
computed maximum solution less than or equal to $1-\epsilon$ is
interpreted to mean that the objective $f(\mathbf{x})$ must be
strictly less than one.  Such qualitative solutions are sufficient for
Aristotelian deduction; for this application the precise value of
$\epsilon$ is not important.  Anyway most optimization solvers use
floating-point arithmetic and rely on various small constants to
control their operation.

To illustrate, consider the problem defined by the premises $BeA$ and
$BiC$.  Referring to Tables \ref{tbl:composites} and
\ref{tbl:descriptions}, the first premise $BeA$ says `$B$ belongs to
no $A$, or there are no $A$' and translates to the probability
constraint $\prob{A,B}=0$ (from which it follows that the conditional
probability $\condprob{B}{A}$ must have the value $0$ unless it is
undefined).  As shown in Table~\ref{tbl:basic-output}~(a) the
probability expression $\prob{A=\tval,B=\tval}$ has the algebraic
value
\begin{math}
x_{1} + x_{2}
\end{math}.  Hence this first premise $BeA$ gives the linear
constraint:
\begin{eqnarray}
x_{1} + x_{2}
& = & 0
\end{eqnarray}
The second premise $BiC$ says `$B$ belongs to some $C$' and translates
to the probability constraint $\prob{B,C}>0$ (requiring that the
conditional probability $\condprob{B}{C}>0$ as well).  As
Table~\ref{tbl:basic-output}~(b) shows,
\begin{math}
  \prob{B=\tval,C=\tval} \Rightarrow
x_{1} + x_{5}
\end{math}.  Thus this second premise $BiC$ asserts the linear
constraint:
\begin{eqnarray}
x_{1} + x_{5}
& > & 0
\end{eqnarray}
Applying epsilon-inequality reformulation gives the modified
constraint
\begin{math}
x_{1} + x_{5}
\ge \epsilon
\end{math}.  Let us choose as our objective function the probability
$\prob{C,\overline{A}}$ that $C$ is true and $A$ is false.  As shown
in Table~\ref{tbl:basic-output}~(c), this probability expression
\begin{math}
  \prob{C=\tval,A=\fval}
\end{math}
evaluates to the algebraic formula
\begin{math}
x_{5} + x_{7}
\end{math}.  

Combining these results, appending the general constraints $0 \le x_i
\le 1$ and $\sum_i x_i = 1$ from the basic probability model in
Section~\ref{sec:basic-model}, and choosing
\begin{math}
  \epsilon =
0.01
\end{math}
to encode strict inequality leads to the following linear optimization
problem to find the minimum feasible value of the chosen objective
probability $\prob{C,\overline{A}}$ subject to the constraints
translated from the categorical premises $BeA$ and $BiC$:
\begin{equation}
\begin{array}{r@{\quad}l}
\mbox{Minimize}: & x_{5} + x_{7} \\
\mbox{subject to}:
& x_{1} + x_{2} + x_{3} + x_{4} + x_{5} + x_{6} + x_{7} + x_{8} = 1 \\
& x_{1} + x_{2} = 0 \\
& x_{1} + x_{5} \ge \epsilon \\
\mbox{and}:
& \epsilon = 0.01 \\
& 0 \le x_{1} \le 1 \\
& 0 \le x_{2} \le 1 \\
& 0 \le x_{3} \le 1 \\
& 0 \le x_{4} \le 1 \\
& 0 \le x_{5} \le 1 \\
& 0 \le x_{6} \le 1 \\
& 0 \le x_{7} \le 1 \\
& 0 \le x_{8} \le 1
\end{array}
  \label{eq:pp1min}
\end{equation}
Linear optimization yields the numerical solution
\begin{math}
0.01
\end{math}.
This computed minimum value says directly that the given constraints
require
\begin{math}
x_{5} + x_{7}
\ge
0.01
\end{math}.
Reversing epsilon-inequality
encoding, we interpret this numerical solution as the strict
algebraic inequality
\begin{math}
x_{5} + x_{7}
  > 0
\end{math} 
which says that the value of the objective function is constrained to
be greater than zero at all feasible points.  The complementary
problem to find the \emph{maximum} value of the same objective subject
to the same constraints yields the solution exactly
\begin{math}
1
\end{math}.
In other words the given constraints require
\begin{math}
x_{5} + x_{7}
\le
1
\end{math}
(which we already knew from the general constraints included in the
basic probability model).
%%   Anyway these results from
%% parametric probability analysis and linear programming can be used to
%% compute Aristotelian syllogism, in a manner described in the following
%% sections.

\subsection{Criteria for Classical and Complementary Syllogism}
\label{sec:criteria}

All four Aristotelian figures use the same query, in which the major
term $A$ is the predicate and the minor term $C$ is the subject.  The
meaning of this common query $A\mathfrak{s}C$, which is shown in
Table~\ref{tbl:figures}, is as follows.  We seek every relation type
$\mathfrak{s}$ from the set
\begin{math}
  \{ a, \acute{a}, e, \acute{e}, i, o, u \}
\end{math}
of composite types from Table~\ref{tbl:composites} for which the
categorical statement $A\mathfrak{s}C$ is a necessary consequence of
the two categorical statements which have been asserted as premises.
There may be zero, one, or many such satisfactory relation types
$\mathfrak{s}$.  If there is at least one satisfactory relation type,
we say that a `syllogism' (deduction) is present; but if there are no
satisfactory relation types then there is no syllogism.  Let us say
that a `classical' syllogism concerns a deduced statement in which the
minor term $C$ is held to be true.  We shall also consider
`complementary' syllogisms in which the minor term $C$ is held to be
false.  Hence to find a complementary syllogism, the query uses the
form $A\mathfrak{s}\overline{C}$ whose subject is the negation of the
minor term.  It is possible to have complementary syllogism with or
without classical syllogism (and likewise to have classical syllogism
with or without complementary syllogism).

\begin{table}
  \sf
  \renewcommand{\arraystretch}{1.25}
  \begin{tabular}{l|c|c|c} \hline
    \renewcommand{\arraystretch}{1}
    \begin{tabular}{@{}l@{}} 
      \thd{Objective} \\ \thd{Probability} 
    \end{tabular} &
    \renewcommand{\arraystretch}{1}
    \begin{tabular}{@{}l@{}} \thd{Algebraic} \\ \thd{Formula} \end{tabular} &
    \renewcommand{\arraystretch}{1}
    \begin{tabular}{@{}l@{}} \thd{Computed} \\ \thd{Minimum} \end{tabular} &
    \renewcommand{\arraystretch}{1}
    \begin{tabular}{@{}l@{}} \thd{Computed} \\ \thd{Maximum} \end{tabular}
    \\ \hline\hline
    \prob{C,A} &
    \begin{math}
x_{1} + x_{3}
    \end{math}
    & $\alpha_1$ & $\beta_1$
    \\ \hline
    \prob{C,\overline{A}} &
    \begin{math}
x_{5} + x_{7}
    \end{math}
    & $\alpha_2$ & $\beta_2$
    \\ \hline
    \prob{\overline{C},A} &
    \begin{math}
x_{2} + x_{4}
    \end{math}
    & $\alpha_3$ & $\beta_3$
    \\ \hline
    \prob{\overline{C},\overline{A}} &
    \begin{math}
x_{6} + x_{8}
    \end{math}
    & $\alpha_4$ & $\beta_4$
    \\ \hline
  \end{tabular}
  \renewcommand{\arraystretch}{1}
  \caption{Notation for minimum and maximum bounds computed on various
    objective functions, for the determination of syllogism.  Each
    bound $\alpha_j$ or $\beta_j$ is the solution of a linear
    optimization problem.}
  \label{tbl:objectives}
\end{table}

The results of certain probability-optimization problems indicate
precisely which relation types are necessary consequences of the
premises provided.  To wit, there are four objective functions whose
minimum and maximum feasible values must be computed, subject to the
constraints translated from the provided premises (and subject also to
the general constraints in the basic probability model that reflect
the laws of probability).  These four objective functions are the
joint probabilities of the various combinations of truth and falsity
of the extreme terms $A$ and $C$:
\begin{displaymath}
    \prob{C,A}
    \qquad
    \prob{C,\overline{A}}
    \qquad
    \prob{\overline{C},A}
    \qquad
    \prob{\overline{C},\overline{A}}
\end{displaymath}
As shown in Table~\ref{tbl:objectives}, the symbols $\alpha_j$ and
$\beta_j$ are used to represent the computed minimum and maximum
values for these four objectives.  Table~\ref{tbl:basic-output}
part~(c) gives the algebraic forms of the relevant probability
expressions, which are also included in Table~\ref{tbl:objectives}.
These algebraic formulas are used as objective functions during the
formulation of optimization problems.  Because these objectives
represent probabilities, each pair of computed minimum and maximum
values is already constrained by $0 \le \alpha_j \le \beta_j \le 1$.
Therefore it would be trivial to compute a minimum value $\alpha_j=0$
or a maximum value $\beta_j=1$ for any queried probability, as the
laws of probability already require these bounds.  A proper deduction
requires the computation of a nontrivial upper or lower bound on at
least one of the queried probabilities.  The numerical solutions
computed depend on the value $\epsilon$ chosen to encode strict
inequality (following the epsilon-inequality encoding scheme discussed
in Section~\ref{sec:linear}), although subsequent interpretation does
not depend on the precise value of $\epsilon$.

The computed minimum and maximum values $\alpha_j$ and $\beta_j$
defined in Table~\ref{tbl:objectives} determine the presence or
absence of syllogism, according to the rules discussed presently.
Table~\ref{tbl:classical-syllogism} shows the criteria for classical
syllogism, describing the necessary relationships between the major
term $A$ as predicate and the (affirmative) minor term $C$ as subject.
Table~\ref{tbl:complementary-syllogism} shows the criteria for
complementary syllogism, describing the necessary relationships
between the major term $A$ as predicate and the negation
$\overline{C}$ of the minor term as subject.  As may be evident
already, the criteria presented in Tables
\ref{tbl:classical-syllogism} and \ref{tbl:complementary-syllogism}
for \emph{deducing} categorical relationships based on relations
involving probabilities are none other than the criteria presented in
Table~\ref{tbl:composites} for \emph{defining} categorical
relationships based on relations involving probabilities, instantiated
for the predicate term $A$ and the subject term $C$ (or its negation
$\overline{C}$).

\renewcommand{\arraystretch}{1.25}

\begin{table}
  \sf
  \begin{tabular}{c|c||l|l} \hline
    \multicolumn{2}{c|}{\bfseries Computed Bounds} & & 
    \bfseries{Categorical}
    \\ \cline{1-2}
    $\prob{C,A}$ & $\prob{C,\overline{A}}$ &
    \bfseries Inferred Probability Relations &
    \bfseries Deduction \\ \hline\hline
    
    & $\beta_2=0$ & 
    $\prob{C,A}=\prob{C}$ & $AaC$ \\ \hline
    
    $\alpha_1>0$ & $\beta_2=0$ & 
    $\prob{C,A}=\prob{C}$ and $\prob{C}>0$ & $A\acute{a}C$ \\ \hline
    
    $\beta_1=0$ & & 
    $\prob{C,A}=0$ & $AeC$ \\ \hline
    
    $\beta_1=0$ & $\alpha_2>0$ & 
    $\prob{C,A}=0$ and $\prob{C}>0$ & $A\acute{e}C$ \\ \hline
    
    $\alpha_1>0$ & & 
    $\prob{C,A} > 0$ & $AiC$ \\ \hline
    
    & $\alpha_2>0$ &
    $\prob{C,A} < \prob{C}$ & $AoC$ \\ \hline
    
    $\alpha_1>0$ & $\alpha_2>0$ &
    $\prob{C,A} > 0$ and $\prob{C,A} < \prob{C}$ & $AuC$ \\ \hline
    
  \end{tabular}
  \caption{Criteria for classical syllogism, describing how the major
    term $A$ must be predicated upon the minor term $C$ when $C$ is
    true, using the given premises as constraints.  Here $\alpha_1$ is
    the minimum and $\beta_1$ the maximum feasible value of
    $\prob{C,A}$; likewise $\alpha_2$ is the minimum and $\beta_2$ the
    maximum feasible value of $\prob{C,\overline{A}}$.  Multiple
    criteria may apply.}
  \label{tbl:classical-syllogism}
\end{table}

\begin{table}
  \sf
  \begin{tabular}{c|c||l|l} \hline
    \multicolumn{2}{c|}{\bfseries Computed Bounds} & & 
    \bfseries{Categorical}
    \\ \cline{1-2}
    $\prob{\overline{C},A}$ & $\prob{\overline{C},\overline{A}}$ &
    \bfseries Inferred Probability Relations &
    \bfseries Deduction \\ \hline\hline
    
    & $\beta_4=0$ & 
    $\prob{\overline{C},A}=\prob{\overline{C}}$ & $Aa\overline{C}$ 
    \\ \hline
    
    $\alpha_3>0$ & $\beta_4=0$ & 
    $\prob{\overline{C},A}=\prob{\overline{C}}$ and
    $\prob{\overline{C}}>0$ &
    $A\acute{a}\overline{C}$ \\ \hline
    
    $\beta_3=0$ & & 
    $\prob{\overline{C},A}=0$ & $Ae\overline{C}$ \\ \hline
    
    $\beta_3=0$ & $\alpha_4>0$ & 
    $\prob{\overline{C},A}=0$ and $\prob{\overline{C}}>0$ & 
    $A\acute{e}\overline{C}$ \\ \hline
    
    $\alpha_3>0$ & & 
    $\prob{\overline{C},A} > 0$ & $Ai\overline{C}$ \\ \hline
    
    & $\alpha_4>0$ &
    $\prob{\overline{C},A} < \prob{\overline{C}}$ & $Ao\overline{C}$ 
    \\ \hline
    
    $\alpha_3>0$ & $\alpha_4>0$ &
    $\prob{\overline{C},A} > 0$ and 
    $\prob{\overline{C},A} < \prob{\overline{C}}$ & 
    $Au\overline{C}$ \\ \hline
    
  \end{tabular}
  \caption{Criteria for complementary syllogism, describing how the
    major term $A$ must be predicated upon the minor term $C$ when $C$
    is false, using the given premises as constraints.  Here
    $\alpha_3$ is the minimum and $\beta_3$ the maximum feasible value
    of $\prob{\overline{C},A}$; likewise $\alpha_4$ is the minimum and
    $\beta_4$ the maximum feasible value of
    $\prob{\overline{C},\overline{A}}$.  Multiple criteria may apply.}
  \label{tbl:complementary-syllogism}
\end{table}

\renewcommand{\arraystretch}{1}

To illustrate, let us consider how the computed optimization result
$\alpha_2>0$ gives the categorical deduction $AoC$, as shown in the
penultimate row of Table~\ref{tbl:classical-syllogism}.  As shown in
Table~\ref{tbl:objectives}, the symbol $\alpha_2$ designates the
minimum feasible value of the probability $\prob{C,\overline{A}}$
subject to the provided constraints.  A computed solution $\alpha_2$
which is strictly greater than zero means that the constraints require
$\prob{C,\overline{A}}>0$.  The laws of probability provide that:
\begin{eqnarray}
  \prob{C} & = & \prob{C,A} + \prob{C,\overline{A}} 
  \label{eq:ca-classical}
\end{eqnarray}
As all probabilities are nonnegative, the inequality
$\prob{C,\overline{A}}>0$ joined with Equation~\ref{eq:ca-classical}
requires in turn that $\prob{C,A}<\prob{C}$.  Returning to
Table~\ref{tbl:composites}, this derived inequality
$\prob{C,A}<\prob{C}$ is precisely the definition of the
\rname{Particular-negative} categorical relationship with predicate
$A$ and subject $C$, abbreviated $AoC$.  Also, as
Table~\ref{tbl:descriptions} shows, this unconditioned-probability
inequality $\prob{C,A}<\prob{C}$ implies the conditional-probability
inequality $\condprob{A}{C}<1$.

Applying this criterion to the example problem from
Section~\ref{sec:linear} produces the deduction that $AoC$ is a
necessary consequence of the premises $BeA$ and $BiC$.
The solutions calculated for the example in
Section~\ref{sec:linear} are the bounds
\begin{math}
  \alpha_2 =
0.01
\end{math}
and
\begin{math}
  \beta_2 =
1
\end{math}
on the objective 
\begin{math}
x_{5} + x_{7}
\end{math}
(translated from the probability expression $\prob{C,\overline{A}}$)
when this objective is subjected to the constraints
\begin{math}
x_{1} + x_{2}
  = 0
\end{math} and
\begin{math}
x_{1} + x_{5}
  > 0
\end{math}
(translated from the categorical premises $BeA$ and $BiC$) and the
constraints $0 \le x_i \le 1$ and $\sum_i x_i = 1$ (from the basic
probability model).  The computed solution
\begin{math}
  \alpha_2 =
0.01
\end{math}
says directly that the given constraints require that the objective
\begin{math}
  \prob{C,\overline{A}} \ge
0.01
\end{math}
at all feasible points.  Following the chain of reasoning outlined
above, we interpret this numerical solution to mean the strict
inequality $\prob{C,\overline{A}} > 0$ and in turn
$\prob{C,A}<\prob{C}$.  This last inequality is the definition of the
categorical relationship $AoC$ (from Table~\ref{tbl:composites}).
Thus is it derived using probability and optimization that the
premises $BeA$ and $BiC$ require the conclusion $AoC$.  In other
words, the premises `$B$ belongs to no $A$, or there are no $A$' and
`$B$ belongs to some $C$' require the conclusion that `The negation of
$A$ belongs to some $C$' (which might alternatively be stated as `$A$
does not belong to some $C$').

Following convention we abbreviate each classical syllogism as
$\mathfrak{m}\mathfrak{n}\mathfrak{s}$-$k$ where $\mathfrak{m}$ is the
code of the major premise's relation type, $\mathfrak{n}$ is the code
for the minor premise, $\mathfrak{s}$ is the code for the deduced
relation type, and $k$ is the number of the problem's figure.  Hence
`aaa-1', `eio-2', and so on.  By analogy each complementary syllogism
is abbreviated as
\begin{math}
  \mathfrak{m}\mathfrak{n}\mathfrak{s}$-$\overline{k}
\end{math}, 
using the bar over the figure number to indicate the negation of the
minor term in the deduced statement.  Syllogistic deductions can also
be displayed in alternative notation using the turnstile symbol or a
tabular arrangement of formulas, as in:
\begin{eqnarray}
  BeA, BiC & \vdash & AoC
\end{eqnarray}
or:
\begin{equation}
  \begin{array}{cl}
    & BeA \\
    & BiC \\ \hline
    \therefore & AoC
  \end{array}
\end{equation}
for the pattern eio-2.

\section{Exhaustive Analysis of Aristotelian Problems}
\label{sec:exhaustive}

Using the probability and optimization methods presented above, let us
now analyze all possible Aristotelian problems with the structure
given in Section~\ref{sec:review} and the types of categorical
relationships enumerated in Table~\ref{tbl:composites}.  For this task
we shall consider all $4$ figures, and within each figure all $7$
possible relationship types for each of the $2$ premises.  This gives
$4 \times 7^2$ or $196$ distinct Aristotelian problems.  For each
Aristotelian problem we set up $8$ optimization problems: one problem
to find the minimum feasible value $\alpha_j$ and one to find the
maximum feasible value $\beta_j$ of each of the $4$ objective
probabilities from Table~\ref{tbl:objectives}.  This gives $8$
optimization problems for each of the $196$ Aristotelian problems,
hence $1,568$ optimization problems altogether.  For each Aristotelian
problem, we compare the results of its $8$ optimization problems with
the deductive criteria listed in Tables \ref{tbl:classical-syllogism}
and \ref{tbl:complementary-syllogism} for classical and complementary
syllogism.  Using each table of deductive criteria, we record a yes/no
answer for whether the optimization results satisfy the requirements
for each of the $7$ types of categorical relationships; thus there are
$14$ yes/no answers for each of the $196$ Aristotelian problems.

The results of this exhaustive analysis are displayed in Tables
\ref{tbl:first-figure} through \ref{tbl:fourth-figure}, in the
following format.  Each result table has two parts: part~(a) which
shows classical syllogisms, and part~(b) which shows complementary
syllogisms.  Within each part there is a $7 \times 7$ inner table.
Each row of an inner table indicates the type of the major premise;
each column indicates the type of the minor premise.  Tabulated within
each cell of an inner table are the codes for all of the valid
deductions from the premises indexed by that row and column, using the
figure indicated.  These valid deductions are syllogisms.

\begin{table}
  \sf
  \begin{tabular}{c@{\quad}|@{\quad}c}
    
    \begin{tabular}{lc}
      & Minor premise ($B\mathfrak{n}C$) \\ 
      \begin{tabular}{@{}l@{}}
        Major \\ premise \\ ($A\mathfrak{m}B$)
      \end{tabular} &
      \begin{tabular}{|c||c|c|c|c|c|c|c|} \hline
        & \bfseries a & 
        \bfseries \'a & 
        \bfseries e & 
        \bfseries \'e & 
        \bfseries i & 
        \bfseries o & 
        \bfseries u \\ \hline\hline
        \bfseries a & a & \'a, a, i & & & i & & i \\ \hline
        \bfseries \'a & a & \'a, a, i & & & i & & i \\ \hline
        \bfseries e & e & \'e, e, o & & & o & & o \\ \hline
        \bfseries \'e & e & \'e, e, o & & & o & & o \\ \hline
        \bfseries i & & & & & & & \\ \hline
        \bfseries o & & & & & & & \\ \hline
        \bfseries u & & & & & & & \\ \hline
      \end{tabular}
    \end{tabular}
    
    &
    
    \begin{tabular}{lc}
      & Minor premise ($B\mathfrak{n}C$) \\ 
      \begin{tabular}{@{}l@{}}
        Major \\ premise \\ ($A\mathfrak{m}B$)
      \end{tabular} &
      \begin{tabular}{|c||c|c|c|c|c|c|c|} \hline
        & \bfseries a & 
        \bfseries \'a & 
        \bfseries e & 
        \bfseries \'e & 
        \bfseries i & 
        \bfseries o & 
        \bfseries u \\ \hline\hline
        \bfseries a & & & & & & & \\ \hline
        \bfseries \'a & & & i & i & & & \\ \hline
        \bfseries e & & & & & & & \\ \hline
        \bfseries \'e & & & o & o & & & \\ \hline
        \bfseries i & & & i & i & & & \\ \hline
        \bfseries o & & & o & o & & & \\ \hline
        \bfseries u & & & u & u & & & \\ \hline
      \end{tabular}
    \end{tabular}
    
    \\ \\
    (a) Classical syllogism ($A\mathfrak{s}C$) & 
    (b) Complementary syllogism ($A\mathfrak{s}\overline{C}$)
  \end{tabular}
  \caption{Analysis of Aristotelian problems in the First Figure.
    Tabulated is every type $\mathfrak{s}$ for which the categorical
    statement $A\mathfrak{s}C$ (classical) or
    $A\mathfrak{s}\overline{C}$ (complementary) is a valid deduction
    from the indexing premises $A\mathfrak{m}B$ and
    $B\mathfrak{n}C$.}
  \label{tbl:first-figure}
\end{table}

For concreteness let us focus on Table~\ref{tbl:first-figure} part~(a)
which describes the classical syllogisms in Aristotle's first figure,
whose major premise $A\mathfrak{m}B$ and minor premise
$B\mathfrak{n}C$ have types $\mathfrak{m}$ and $\mathfrak{n}$.
Focusing more specifically on one particular problem in the first
figure, the third row (labeled `$e$') indicates the major premise
$AeB$, the first column (labeled `$a$') indicates the minor premise
$BaC$.  The solitary cell entry `$e$' at this row and column says that
the only valid deduction from these premises is $AeC$.  That is, from
the premises `$A$ belongs to no $B$, or there are no $B$' and `$B$
belongs to all $C$, or there are no $C$,' there follows the conclusion
`$A$ belongs to no $C$, or there are no $C$.'  This syllogistic
pattern `eae-1' is also known by the medieval name `Celarent'.  (The
vowels in the medieval names give codes for categorical relations in
the same sequence as the $\mathfrak{m}\mathfrak{n}\mathfrak{s}$-$k$
abbreviation.  The consonants in the medieval names also convey
information, which is not essential to the analysis presented here.)

Remaining in the third row but moving over to the second column
(labeled `$\acute{a}$') gives slightly different result.  Here there
are three cell entries $\acute{e}$, $e$, and $o$.  These indicate that
all three statements $A\acute{e}C$, $AeC$, and $AoC$ are valid
deductions from the premises $AeB$ and $B\acute{a}C$.  That is, from
the premises `$A$ belongs to no $B$, or there are no $B$' and
`$B$ belongs to all $C$, and there are some $C$,' there follow three
necessary conclusions:
\begin{description}
  \item[$A\acute{e}C$:] `$A$ belongs to no $C$, and there are some $C$'
  \item[$AeC$:] `Either $A$ belongs to no $C$, or there are no $C$'
  \item[$AoC$:] `The negation of $A$ belongs to some $C$'
\end{description}
Let us abbreviate these patterns of deduction as `e\'a\'e-1',
`e\'ae-1', and `e\'ao-1', with the accent marks used as in
Table~\ref{tbl:composites} to distinguish \rname{existential} from
\rname{material} statements.  By analogy with the medieval names we
could call these patterns of syllogism `Cel\'ar\'ent', `Cel\'arent'
and `Cel\'aront'.  Note that there is no syllogism `eao-1' listed in
Table~\ref{tbl:first-figure}~(a); the pattern named `Celaront' would
be a mistake.  In interpreting these results it is important to pay
attention to what is absent as well as to what is present.

\begin{table}
  \sf
  \begin{tabular}{c@{\quad}|@{\quad}c}
    
    \begin{tabular}{lc}
      & Minor premise ($B\mathfrak{n}C$) \\ 
      \begin{tabular}{@{}l@{}}
        Major \\ premise \\ ($B\mathfrak{m}A$)
      \end{tabular} &
      \begin{tabular}{|c||c|c|c|c|c|c|c|} \hline
        & \bfseries a & 
        \bfseries \'a & 
        \bfseries e & 
        \bfseries \'e & 
        \bfseries i & 
        \bfseries o & 
        \bfseries u \\ \hline\hline
        \bfseries a & & & e & \'e, e, o & & o & o \\ \hline
        \bfseries \'a & & & e & \'e, e, o & & o & o \\ \hline
        \bfseries e & e & \'e, e, o & & & o & & o \\ \hline
        \bfseries \'e & e & \'e, e, o & & & o & & o \\ \hline
        \bfseries i & & & & & & & \\ \hline
        \bfseries o & & & & & & & \\ \hline
        \bfseries u & & & & & & & \\ \hline
      \end{tabular}
    \end{tabular}
    
    &
    
    \begin{tabular}{lc}
      & Minor premise ($B\mathfrak{n}C$) \\ 
      \begin{tabular}{@{}l@{}}
        Major \\ premise \\ ($B\mathfrak{m}A$)
      \end{tabular} &
      \begin{tabular}{|c||c|c|c|c|c|c|c|} \hline
        & \bfseries a & 
        \bfseries \'a & 
        \bfseries e & 
        \bfseries \'e & 
        \bfseries i & 
        \bfseries o & 
        \bfseries u \\ \hline\hline
        \bfseries a & & & & & & & \\ \hline
        \bfseries \'a & & & i & i & & & \\ \hline
        \bfseries e & & & & & & & \\ \hline
        \bfseries \'e & i & i & & & & & \\ \hline
        \bfseries i & & & i & i & & & \\ \hline
        \bfseries o & i & i & & & & & \\ \hline
        \bfseries u & i & i & i & i & & & \\ \hline
      \end{tabular}
    \end{tabular}
    
    \\ \\
    (a) Classical syllogism ($A\mathfrak{s}C$) & 
    (b) Complementary syllogism ($A\mathfrak{s}\overline{C}$)
  \end{tabular}
  \caption{Analysis of Aristotelian problems in the Second Figure.
    Tabulated is every type $\mathfrak{s}$ for which the categorical
    statement $A\mathfrak{s}C$ (classical) or
    $A\mathfrak{s}\overline{C}$ (complementary) is a valid deduction
    from the indexing premises $B\mathfrak{m}A$ and $B\mathfrak{n}C$.}
  \label{tbl:second-figure}
\end{table}

\begin{table}
  \sf
  \begin{tabular}{c@{\quad}|@{\quad}c}
    
    \begin{tabular}{lc}
      & Minor premise ($C\mathfrak{n}B$) \\
      \begin{tabular}{@{}l@{}}
        Major \\ premise \\ ($A\mathfrak{m}B$)
      \end{tabular} &
      \begin{tabular}{|c||c|c|c|c|c|c|c|} \hline
        & \bfseries a & 
        \bfseries \'a & 
        \bfseries e & 
        \bfseries \'e & 
        \bfseries i & 
        \bfseries o & 
        \bfseries u \\ \hline\hline
        \bfseries a &  & i & & & i & & i \\ \hline
        \bfseries \'a & i & i & & & i & & i \\ \hline
        \bfseries e & & o & & & o & & o \\ \hline
        \bfseries \'e & o & o & & & o & & o \\ \hline
        \bfseries i & i & i & & & & & \\ \hline
        \bfseries o & o & o & & & & & \\ \hline
        \bfseries u & u & u & & & & & \\ \hline
      \end{tabular}
    \end{tabular}
    
    &
    
    \begin{tabular}{lc}
      & Minor premise ($C\mathfrak{n}B$) \\
      \begin{tabular}{@{}l@{}}
        Major \\ premise \\ ($A\mathfrak{m}B$)
      \end{tabular} &
      \begin{tabular}{|c||c|c|c|c|c|c|c|} \hline
        & \bfseries a & 
        \bfseries \'a & 
        \bfseries e & 
        \bfseries \'e & 
        \bfseries i & 
        \bfseries o & 
        \bfseries u \\ \hline\hline
        \bfseries a & & & & i & & i & i \\ \hline
        \bfseries \'a & & & i & i & & i & i \\ \hline
        \bfseries e & & & & o & & o & o \\ \hline
        \bfseries \'e & & & o & o & & o & o \\ \hline
        \bfseries i & & & i & i & & & \\ \hline
        \bfseries o & & & o & o & & & \\ \hline
        \bfseries u & & & u & u & & & \\ \hline
      \end{tabular}
    \end{tabular}
    
    \\ \\
    (a) Classical syllogism ($A\mathfrak{s}C$) & 
    (b) Complementary syllogism ($A\mathfrak{s}\overline{C}$)
  \end{tabular}
  \caption{Analysis of Aristotelian problems in the Third
    Figure.  Tabulated is every type $\mathfrak{s}$ for
    which the categorical statement $A\mathfrak{s}C$ (classical) or
    $A\mathfrak{s}\overline{C}$ (complementary) is a valid deduction
    from the indexing premises $A\mathfrak{m}B$ and
    $C\mathfrak{n}B$.}
  \label{tbl:third-figure}
\end{table}

\begin{table}
  \sf
  \begin{tabular}{c@{\quad}|@{\quad}c}
    
    \begin{tabular}{lc}
      & Minor premise ($C\mathfrak{n}B$) \\
      \begin{tabular}{@{}l@{}}
        Major \\ premise \\ ($B\mathfrak{m}A$)
      \end{tabular} &
      \begin{tabular}{|c||c|c|c|c|c|c|c|} \hline
        & \bfseries a & 
        \bfseries \'a & 
        \bfseries e & 
        \bfseries \'e & 
        \bfseries i & 
        \bfseries o & 
        \bfseries u \\ \hline\hline
        \bfseries a & & & e & e & & & \\ \hline
        \bfseries \'a & i & i & e & e & & & \\ \hline
        \bfseries e & & o & & & o & & o \\ \hline
        \bfseries \'e & & o & & & o & & o \\ \hline
        \bfseries i & i & i & & & & & \\ \hline
        \bfseries o & & & & & & & \\ \hline
        \bfseries u & i & i & & & & & \\ \hline
      \end{tabular}
    \end{tabular}
    
    &
    
    \begin{tabular}{lc}
      & Minor premise ($C\mathfrak{n}B$) \\
      \begin{tabular}{@{}l@{}}
        Major \\ premise \\ ($B\mathfrak{m}A$)
      \end{tabular} &
      \begin{tabular}{|c||c|c|c|c|c|c|c|} \hline
        & \bfseries a & 
        \bfseries \'a & 
        \bfseries e & 
        \bfseries \'e & 
        \bfseries i & 
        \bfseries o & 
        \bfseries u \\ \hline\hline
        \bfseries a & e & e & & & & & \\ \hline
        \bfseries \'a & e & e & i & i & & & \\ \hline
        \bfseries e & & & & o & & o & o \\ \hline
        \bfseries \'e & & & & o & & o & o \\ \hline
        \bfseries i & & & i & i & & & \\ \hline
        \bfseries o & & & & & & & \\ \hline
        \bfseries u & & & i & i & & & \\ \hline
      \end{tabular}
    \end{tabular}
    
    \\ \\
    (a) Classical syllogism ($A\mathfrak{s}C$) & 
    (b) Complementary syllogism ($A\mathfrak{s}\overline{C}$)
  \end{tabular}
  \caption{Analysis of Aristotelian problems in the Fourth
    Figure.  Tabulated is every type $\mathfrak{s}$ for
    which the categorical statement $A\mathfrak{s}C$ (classical) or
    $A\mathfrak{s}\overline{C}$ (complementary) is a valid deduction
    from the indexing premises $B\mathfrak{m}A$ and
    $C\mathfrak{n}B$.}
  \label{tbl:fourth-figure}
\end{table}

\subsection{Classical Modes of Syllogism Reproduced}

The results of probability-optimization analysis shown in Tables
\ref{tbl:first-figure} through \ref{tbl:fourth-figure} reproduce the
standard modes of Aristotelian syllogism, which are described for
example in \cite{aristotle-prior-smith} and
\cite{sep-medieval-syllogism}.  Let us focus for a moment on the four
relation types $\{a,e,i,o\}$, excluding the
\rname{particular-intermediate} relation $u$ and also excluding the
\rname{existential} subtypes $\acute{a}$ and $\acute{e}$ of the
\rname{universal-affirmative} and \rname{universal-negative}
relations.  Table~\ref{tbl:first-figure}~(a) indicates the following
four patterns of syllogism for Aristotle's first figure (displayed
here with their medieval names):
\begin{quote}
  aaa-1 (Barbara) \quad
  aii-1 (Darii) \quad
  eae-1 (Celarent) \quad
  eio-1 (Ferio)
\end{quote}
Remember, for each syllogism denoted
$\mathfrak{m}\mathfrak{n}\mathfrak{s}$-$k$, the type $\mathfrak{m}$ of
the major premise gives the row, the type $\mathfrak{n}$ of the minor
premise gives the column, and the type $\mathfrak{s}$ of the deduced
statement appears in the cell at that row and column in the specified
results table; $k$ is the number of the figure.
Table~\ref{tbl:second-figure}~(a) shows the following four modes of
syllogism for the second figure:
\begin{quote}
  aee-2 (Camestres) \quad
  aoo-2 (Baroco) \quad
  eae-2 (Cesare) \quad
  eio-2 (Festino)
\end{quote}
Table~\ref{tbl:third-figure}~(a) shows the following four modes of
syllogism for the third figure:
\begin{quote}
  aii-3 (Datisi) \quad
  eio-3 (Ferison) \quad
  iai-3 (Disamis) \quad
  oao-3 (Bocardo)
\end{quote}
Table~\ref{tbl:fourth-figure}~(a) shows the following three modes of
syllogism for the fourth figure:
\begin{quote}
  aee-4 (Camenes) \quad
  eio-4 (Fresison) \quad
  iai-4 (Dimaris)
\end{quote}
There are few notable absences from the modes of syllogism listed
here: aai-3 (Darapti); eao-3 (Felapton); aai-4 (Bramantip); and eao-4
(Fesapo).  These cases are addressed in the next section.

\subsection{Existential Fallacies Revealed}

Discipline about the different subtypes of universal statements
prevents existential fallacies from contaminating our analysis of
Aristotelian problems.  For example, Table~\ref{tbl:third-figure}~(a)
shows that there are no valid syllogisms `aai-3' nor `eao-3'.  In
other words the patterns `Darapti' and `Felapton' are invalid using
the \rname{material} interpretation of their \rname{universal}
premises.  However the patterns `D\'arapti', `Dar\'apti', and
`D\'ar\'apti' are all valid, as are `F\'elapton', `Fel\'apton', and
`F\'el\'apton'.  Valid syllogism in these cases requires that the
truth of the middle term $B$ is not impossible \emph{a priori}, in
other words the existence of $B$ (in the third figure, the middle term
$B$ is the subject of both premises).

Likewise, Table~\ref{tbl:fourth-figure}~(a) shows that there are no
valid syllogisms `aai-4' nor `eao-4'.  In other words the patterns
`Bramantip' and `Fesapo' are invalid using the \rname{material}
interpretation of their \rname{universal} premises.  However the
following patterns are valid: `Br\'amantip', `Br\'am\'antip',
`Fes\'apo', and `F\'es\'apo'.  For the `Bramantip' patterns, syllogism
requires that the truth of the major term $A$ is not impossible
\emph{a priori} (that is, they require the existence of $A$, which is
the predicate of the major premise).  For the `Fesapo' patterns,
syllogism requires that the truth of the middle term $B$ is not
impossible \emph{a priori} (that is, they require the existence of
$B$, which is the predicate of the minor premise).

\subsection{Complementary Syllogisms Added}

As an additional benefit, analysis according to the
probability-optimization paradigm finds instances of complementary
syllogism which were not heretofore appreciated.  Complementary
syllogisms can recover information contained within the premises that
would otherwise be lost.  For example, consider the problem in the
first figure with premises $AiB$ and $BeC$.  There is no classical
syllogism in this case: no particular relation is required when $A$ is
predicated on $C$.  However there is a complementary syllogism in
which the negation of the minor term $C$ is used as the subject of the
query.  As Table~\ref{tbl:first-figure}~(b) shows, the deduction
$Ai\overline{C}$ is a necessary consequence of the premises $AiB$ and
$BeC$.  That is, there follows from the premises `$A$ belongs to some
$B$' and `$B$ belongs to no $C$, or there are no $C$' the necessary
consequence that `$A$ belongs to some non-$C$'.  This instance of
complementary syllogism is abbreviated `iei-\={1}' (note the bar over
the figure number; there is no valid classical syllogism iei-1).  As
it happens, all four patterns iei-\={1}, iei-\={2}, iei-\={3}, and
iei-\={4} represent valid deductions:
\begin{eqnarray}
  AiB, BeC & \vdash & Ai\overline{C} \\
  BiA, BeC & \vdash & Ai\overline{C} \\
  AiB, CeB & \vdash & Ai\overline{C} \\
  BiA, CeB & \vdash & Ai\overline{C}
\end{eqnarray}
Despite their different figures, each pair of premises here shares the
identical (and solitary) consequence that `$A$ belongs to some
non-$C$'.

\section{Conclusion}

Using probability and optimization, it is possible to compute
solutions to the logic problems that Aristotle described in
\emph{Prior Analytics}.  The requisite calculations take advantage of
two kinds of mappings: first between categorical statements and
relations involving probabilities; and second between probability
expressions and algebraic expressions (linear functions and fractional
linear functions).  These mappings allow categorical statements to be
translated to and from linear equalities and inequalities involving a
few real-valued variables.  To begin the analysis, Aristotelian
premises are translated into systems of linear constraints.  Numerical
bounds are then computed on the feasible values of certain objective
functions, subject to these constraints.  These computed bounds reveal
precisely which categorical statements are necessary consequences of
the premises that were asserted.  Every valid syllogism from an
Aristotelian problem can be computed in this way.

There are several benefits to this probability-optimization
formulation of Aristotle's logic.  First, the inference uses quite
ordinary mathematical methods: symbolic probability inference (which
is essentially arithmetic with polynomial expressions) and linear
programming.  It is straightforward to write computer programs to
automate these calculations (as the author has done to generate the
result which are reported above).  Second, the results of
probability-optimization analysis not only reproduce the known modes
of Aristotelian syllogism; they also add new deductive results.  The
analysis clarifies the role of existential import in certain patterns
of syllogism (such as the incorrect modes aai-3, eao-3, aai-4, and
eao-4).  The analysis adds new appreciation for `complementary'
syllogisms: deduced consequences in which the subject is held to be
false instead of true.  For example, it happens that there is no valid
\emph{classical} syllogism from a \rname{Particular-affirmative} major
premise (code $i$) and a \rname{Universal-negative} minor premise
(whether \rname{material} or \rname{existential}, code $e$ or
$\acute{e}$) in any figure.  But in every figure this combination of
premise types leads to a valid \emph{complementary} syllogism: the
major term $A$ must have the \rname{Particular-affirmative} relation
to the negation $\overline{C}$ of the minor term (the shared deduction
$Ai\overline{C}$ says that `$A$ belongs to some non-$C$').

There are many ways to extend this computational framework to provide
even more capabilities.  In addition to computing the categorical
relationships that are \emph{necessary} consequences of the given
premises, it is possible to compute the relationships that are merely
\emph{potential} consequences of the given premises, as well as those
relationships that are \emph{inconsistent} with the premises (using
criteria modified from those given in Tables
\ref{tbl:classical-syllogism} and \ref{tbl:complementary-syllogism}).
Thus alethic modalities of truth can be inferred by numerical
computation.  Furthermore, it is possible to use any numbers of terms
and premises, with the terms distributed among the premises and query
in an arbitrary fashion.  It is not necessary to stick to Aristotle's
original restrictions of using two premises with a common middle term
that is not included in the query.  Finally, it is possible to use
probability models other than the basic one introduced here.  As
discussed further in \cite{norman-conditionals}, different probability
models (with the full-joint probability distribution over categorical
terms factored into several input tables) allow the use of other
semantic types of conditional statements (such as subjunctives), at
the expense of introducing nonlinearity into the polynomials used.

\begin{raggedright}
  \small
  \bibliographystyle{plain}
  \bibliography{../decisions}
%  The author's preprints are available on the Internet at
%  \url{http://www.arXiv.org/a/norman_j_1}.
\end{raggedright}

%% \end{document}
%% %!! exit

\clearpage
\appendix

\section{Source Code}

The probability and optimization results displayed above were
generated by the author's preprocessor \texttt{pqlpp} in response to
commands embedded as comments in a specially prepared \LaTeX{} source
file.  The preprocessor replaced the commands with their output, and
the resulting file was processed by the usual \LaTeX{} tools to make
the document that you see now.  The source code for the probability
model, and instructions for performing the exhaustive analysis
reported in Section~\ref{sec:exhaustive}, are presented in this
appendix.  The results of exhaustive analysis were generated as a CSV
(comma-separated value) text file (viewable as a spreadsheet), whose
contents the author manually copied into the tables in
Section~\ref{sec:exhaustive}.

\subsection{Probability Model: \texttt{aristotle.pql}}

This is the specification of the parametric probability network used
for analysis, in the format used by the author's \texttt{pqlsh} and
\texttt{pqlpp} programs.

\begin{footnotesize}
\begin{verbatim}
// aristotle.pql: setup for Aristotle's syllogism, just add constraints

parameter e { 
  tex = "\epsilon"; 
  label = "Strict inequality"; 
  range = 0.01; 
}

primary A { label = "Major term"; states = binary; }
primary B { label = "Middle term"; states = binary; }
primary C { label = "Minor term"; states = binary; }

clique _C; potential ( _C : A B C ) { parametric(x); }
\end{verbatim}
\end{footnotesize}

\subsection{Script for Exhaustive Analysis: \texttt{syllogism.tcl}}

This is the TCL script for performing exhaustive analysis, using the
author's \texttt{pqlsh} program.

\begin{footnotesize}
\begin{verbatim}
# Analyze all possible Aristotlean syllogism problems, as linear programs

proc premise_constraint { model tbl type epsilon } {

    # e.g. AaB means Pr(B=T,A=F)==0 hence Pr(A|B)==1
    # input $tbl has Pr(B,A) for premise AmB; first 2 items have B=T
    # type `aa' is universal affirmative with existential import
    # AaaB means Pr(B=T,A=F)==0 && Pr(B=T,A=T)>0, hence Pr(A|B)==1 && Pr(B)>0
    # type `ee' is universal negative with existential import
    # AeeB means Pr(B=T,A=T)==0 && Pr(B=T,A=F)>0, hence Pr(A|B)==0 && Pr(B)>0

    switch $type {
	a { $model add_constraint "[$tbl item 2] == 0" }
	aa { $model add_constraint "[$tbl item 2] == 0";
	    $model add_constraint "[$tbl item 1] >= $epsilon" }
	e { $model add_constraint "[$tbl item 1] == 0" }
	ee { $model add_constraint "[$tbl item 1] == 0";
	    $model add_constraint "[$tbl item 2] >= $epsilon" }
	i { $model add_constraint "[$tbl item 1] >= $epsilon" }
	o { $model add_constraint "[$tbl item 2] >= $epsilon" }
	u { $model add_constraint "[$tbl item 1] >= $epsilon";
	    $model add_constraint "[$tbl item 2] >= $epsilon" }
	default { puts "premise_constraint: unknown type $type ignored" }
    }
}

proc conclusions { uval tbl_lb tbl_ub } {
    upvar $uval val

    foreach n {1 2 3 4} { 
	set alpha($n) [$tbl_lb item $n]; set beta($n) [$tbl_ub item $n]
    }

    # first test for classical syllogism
    set val(a) [expr $beta(2) == 0]
    set val(e) [expr $beta(1) == 0]
    set val(i) [expr $alpha(1) > 0]
    set val(o) [expr $alpha(2) > 0]
    set val(u) [expr $alpha(1) > 0 && $alpha(2) > 0]
    # ...with direct existential import for universals
    set val(aa) [expr $beta(2) == 0 && $alpha(1) > 0]
    set val(ee) [expr $beta(1) == 0 && $alpha(2) > 0]

    # next test for complementary syllogism
    set val(ac) [expr $beta(4) == 0]
    set val(ec) [expr $beta(3) == 0]
    set val(ic) [expr $alpha(3) > 0]
    set val(oc) [expr $alpha(4) > 0]
    set val(uc) [expr $alpha(3) > 0 && $alpha(4) > 0]
    # ...with direct existential import for universals
    set val(aac) [expr $beta(4) == 0 && $alpha(3) > 0]
    set val(eec) [expr $beta(3) == 0 && $alpha(4) > 0]

    # build a summary string
    set res ""
    foreach t {aa ee a e u i o} { if $val($t) {lappend res $t} }
    foreach t {aac eec ac ec uc ic oc} { if $val($t) {lappend res $t} }
    return $res
}

proc syllogism { model {A "A"} {B "B"} {C "C"} {epsilon "e"} {sep ","} } {
    set fig1 [dict create major "$A $B" minor "$B $C"]
    set fig2 [dict create major "$B $A" minor "$B $C"]
    set fig3 [dict create major "$A $B" minor "$C $B"]
    set fig4 [dict create major "$B $A" minor "$C $B"]
    set figures [dict create 1 $fig1 2 $fig2 3 $fig3 4 $fig4]

    # adjust this as needed e.g. to ignore aa, ee, u
    set types {a aa e ee i o u}
    
    set tbl_query [[$model table $C $A] infer];	# for bounds on Pr(A|C)

    puts [join "id figure pmajor pminor conseq" $sep]

    dict for {fig terms} $figures {

	# predicates and subjects of major and minor premises
	set major_pred [lindex [dict get $terms major] 0]
	set major_subj [lindex [dict get $terms major] 1]
	set minor_pred [lindex [dict get $terms minor] 0]
	set minor_subj [lindex [dict get $terms minor] 1]

	# joint prob tables: e.g. for AeB, regarding Pr(A|B), infer Pr(B,A)
	set tbl_major [[$model table $major_subj $major_pred] infer]
	set tbl_minor [[$model table $minor_subj $minor_pred] infer]

	set n 0
	foreach major_type $types {
	    foreach minor_type $types {
		incr n; set id [format "%.2f" [expr $fig + $n/100.0]];

		# formluate and solve optimization problems
		$model clear_constraints
		premise_constraint $model $tbl_major $major_type $epsilon
		premise_constraint $model $tbl_minor $minor_type $epsilon
		$model ppinit
		set tbl_lb [$tbl_query lower_bound]
		set tbl_ub [$tbl_query upper_bound]
		set tqb [$tbl_query bound]

		# determine conclusions
		array set cq {}
		set cqsum [conclusions cq $tbl_lb $tbl_ub]

		# print results to standard output
		puts -nonewline [join "$id $fig $major_type $minor_type" $sep]
		puts "$sep$cqsum"

		# print details to standard error
		puts -nonewline stderr "\n\nARISTOTLE $id: ";
		puts -nonewline stderr "$major_pred $major_type $major_subj, ";
		puts stderr "$minor_pred $minor_type $minor_subj, ";
		puts stderr [$model constraints]
		puts stderr [$tqb print]
	    }
	}
   }
}

set m [pql::load aristotle.pql]
pql::print_format text
syllogism $m A B C e ","
\end{verbatim}
\end{footnotesize}

\clearpage
\tableofcontents

\end{document}